# Symmetric Near-Field Schur's Complement Preconditioner for Hierarchal Electric Field Integral Equation Solver

Yoginder Kumar Negi [1], N. Balakrishnan [1], Sadasiva M. Rao [2]

[1] Supercomputer Education Research Centre, Indian Institute of Science, Bangalore, India
[2] Naval Research Laboratory, Washington DC, USA
[*] yknegi@gmail.com

**Abstract:** In this paper, a robust and effective preconditioner for the fast Method of Moments (MoM) based Hierarchal Electric Field Integral Equation (EFIE) solver is proposed using symmetric near-field Schur's complement method. In this preconditioner, near-field blocks are scaled to a diagonal block matrix and these near-field blocks are replaced with the scaled diagonal block matrix which reduces the near-field storage memory and the overall matrix vector product time. Scaled diagonal block matrix is further used as a preconditioner and due to the block diagonal form of the final preconditioner, no additional fill-ins are introduced in its inverse. The symmetric property of the near-field blocks is exploited to reduce the preconditioner setup time. Near linear complexity of preconditioner set up and solve times is achieved by near-field block ordering, using graph bandwidth reduction algorithms and compressing the fill-in blocks in preconditioner computation. Preconditioner set up time is reduced to half by using the symmetric property and near-field block ordering. It has been shown using a complexity analysis that the cost of preconditioner construction in terms of computation and memory is linear. Numerical experiments demonstrate an average of 1.5 - 2.3x speed-up in the iterative solution time over Null-Field based preconditioners.

## 1. Introduction

With the continuous increase in compute power and system memory, the need for solving large and complex electromagnetics problems accurately and efficiently has grown in the last few decades. Few of the complex problems are real world Radar Cross Section (RCS) computation, complex antenna and system level design. System level design technology trends like System-in-Package (SiP) and System-on-Chip (SoC) are becoming increasingly complex and relevant. To ensure the functionality and reliability in chip-to-chip communication and system-level performance, package-board electrical performance parameters like signal integrity (SI), power integrity (PI), electromagnetic interference (EMI) are critical. Obviously, for solving all these complex and large problems an accurate and stable three-dimensional (3D) fast solver is the prime requirement.

The 3D full-wave computational tools have gained in importance over the last couple of decades. Based on the formulation used, 3D electromagnetic solvers can be classified under differential and integral forms. Differential equation based methods include Finite Element Method (FEM)[1] and Finite Difference Time Domain (FDTD) [2] whereas integral equation based techniques include Boundary Element Method (BEM) or Method of Moments (MoM) [3]. FEM and FDTD require a volumetric discretization bounded by an absorbing boundary or a perfectly matched layer (PML). This leads to an increase in matrix size. Integral equation-based methods, on the other hand, require only a surface mesh leading to a smaller number of solution unknowns. However, due to the Green's Function interactions, the MoM system matrix is dense, and with increasing problem size, the solution presents a time-memory bottleneck. As shown in Table 1, for a matrix size of $N \times N$, the conventional LU-based direct solver requires $O(N^3)$ LU factorize time and $O(N^2)$ solve time, whereas a regular iterative solver takes $O(N^2) \times p \times n_p$ for $n_p$ right-hand-sides (RHS) and $p$ is an average number of iterations per RHS. Over the last few decades, several fast-iterative solver techniques have been developed with the complexity of $O(NlogN)$ matrix fill time and matrix vector product time for each iteration. Few of the fast solvers are Multilevel Fast Multipole Algorithm (MLFMA) [4-5], low-rank based methods like SVD/QR compression method [6-7] and Multilevel Adaptive Cross Approximation (MLACA) method [8-9], pre-corrected FFT-based technique [10] and Adaptive Integral Method (AIM) [11]. In this work, we have used multilevel re-compressed Adaptive Cross Approximation (ACA) based Hierarchal Matrix (H-Matrix) Electric Field Integral Equation (EFIE) solver.

**Table 1.** Matrix Setup and Solve Time Complexity

| Method | Setup Time | Solve Time |
| --- | --- | --- |
| Direct | $O(N^2)$ | $O(N^3) + O(N^2) \times n_p$ |
| Iterative | $O(N^2)$ | $O(N^2) \times p \times n_p$ |
| Fast Iterative | $O(NlogN)$ | $O(NlogN) \times p \times n_p$ |

The fast-solver algorithms [4-11] accelerate the MoM matrix vector product to near linear complexity under the framework of a Krylov subspace based iterative solution. Therefore, the solution time for these methods are highly dependent on the number of iterations (*p*) required for convergence. The rate of convergence depends on the system-matrix condition number or more specifically, on the distribution of matrix eigen-values. For the EFIE [12]-based MoM solution, the fast solver iterative solution for large unknown leads to a high number of iterations which results in large solve-times particularly for multiple RHS or multiport excitation chip-package-system analysis.

However, The EFIE system-matrix may become ill-conditioned when solving electrically small structure size [13] leading to a low frequency breakdown, which can be addressed by incorporating loop-star or loop tree basis function [14-15]. Further, high mesh density [16] and poor quality mesh [17] can also lead to ill-conditioning of the MoM matrix, which may be rectified by Calderon preconditioners [18]. In a realistic scenario, all structures are closed and may lead to an ill-conditioned EFIE matrix due to internal resonance [19, 20]. This can be addressed by using Combined Field Integral Equations (CFIE) [19, 20].

In this work, a system-matrix preconditioner that can be used in conjunction with other preconditioning schemes like low-frequency or dense-mesh preconditioners is presented. The proposed method can be used either as a left (1) or right (2) preconditioned system:

$$[P^{-1}][Z]x = [P^{-1}] b \quad (1)$$

$$[Z][P^{-1}]\tilde{x} = b \quad (2)$$

Where $P$ is the preconditioner matrix, $Z$ is the EFIE MoM matrix, $b$ is excitation vector, $x$ and $\tilde{x}$ are the solution vectors, where $\tilde{x}$ is $[P]x$. Different system-matrix preconditioning techniques for fast MoM solution have been proposed in the literature [21-24]. The diagonal preconditioner and block diagonal preconditioner have low setup time but their effectiveness is limited for large scale problems with high iteration count. Approximate Inverse methods like the Sparse Approximate Inverse (SAI) method [22] tries to predict $[P^{-1}]$ directly based on the least square method. They are particularly attractive for parallel applications but is sometimes limited by the high preconditioner setup cost. In a threshold-based Incomplete LU (ILUT) method [23, 24] fill-ins are dropped based on a threshold and/or a fill-in count. The effectiveness of ILUT is limited by the selection of the control parameters, threshold value and fill-in count which may not be uniform for different problem types. Also due to the sparse LU factorization, ILUT is serial in nature, therefore restricting its use in a parallel framework. Typically, an inverse of the near-field MoM blocks can act as a good preconditioner. Block diagonal forms of preconditioner is proposed in [25, 26] where the null-field method is used to scale the symmetric near-field matrix to diagonal blocks. In the case of the null-field method, the fill-in blocks for scaling is not considered which limits the effective scaling of the near-field to diagonal block thus giving high storage and matrix vector product time. Schur's complement [29] is a well-known method in linear algebra for block diagonalization of matrices and has been applied to form a domain decomposition method for symmetric MoM [28]. As preconditioners, the Schur's complement method has been applied in conjunction with SAI [29] and for sparse FEM matrices [30].

In this paper, building upon [31], a Schur's complement preconditioner with near linear time and memory complexity is proposed in conjunction with a re-compressed ACA [32, 33] based fast edge based MoM [12] solution. In the proposed method symmetric near-field is completely block diagonalized by using the Schur's complement technique. The complexity of the preconditioner set up time is reduced by adapting bandwidth reduction graph-based algorithms and exploiting the symmetric property of the Schur's method. Further, compressing the fill-in blocks in the scaling coefficient matrix vector product time and storage can be reduced. The new scaled near-field is block diagonal and can be used as an effective preconditioner since no fill-in is introduced during LU factorization of the preconditioner matrix. Since Schur's complement is known to be efficient for parallelization [34], the preconditioner setup process can be parallelized efficiently. The proposed method is purely algebraic and can be applied to higher order MoM like Locally Corrected Nyström (LCN) method implemented with EFIE, MFIE or CFIE [35, 36].

The paper is structured as follows: in section 2, a brief description of 3D full-wave MoM re-compressed ACA based H-Matrix is presented. In section 3, the proposed near-field Schur's preconditioner is presented along with the method to reduce the preconditioner set up time and matrix vector product time. In section 4 complexity analysis of computation cost and memory requirement for preconditioner, storage is presented and in section 5, the efficiency and accuracy of the proposed preconditioner are validated by comparing the results with null-field preconditioner, ILUT and a commercial solver. Section 6 concludes the paper.

## 2. Re-Compressed ACA based Hierarchal 3D MoM Solver

An arbitrary shaped 3D conductor electromagnetic problem can be solved using full-wave EFIE based MoM. In EFIE boundary condition on the conducting surface ($s$) is given as:

$$[\vec{E}_s + \vec{E}_i]_{tan} = 0 \quad (3)$$

Where $\vec{E}_s$ is the scattered electric field produced by the currents induced on the conducting body, $\vec{E}_i$ is the incident electric field and the subscript $tan$ denoted the tangential component on the scattering surface. Scattered electric field can be written as

$$\vec{E}_s = -j\omega \vec{A}(\vec{r}) - \nabla \phi(\vec{r}) \quad (4)$$

Where $\omega$ is angular frequency, $\vec{A}(\vec{r})$ is vector potential and $\phi(\vec{r})$ is scalar potential. Potentials $\vec{A}(\vec{r})$ and $\phi(\vec{r})$ are given as

$$\vec{A}(\vec{r}) = \frac{\mu}{4\pi} \int \frac{e^{-jk|\vec{r}-\vec{r}'|}}{|\vec{r}-\vec{r}'|} \vec{J}(\vec{r}')ds' \quad (5)$$

$$\phi(\vec{r}) = \frac{1}{4\pi\varepsilon} \int \frac{e^{-jk|\vec{r}-\vec{r}'|}}{|\vec{r}-\vec{r}'|} \rho(\vec{r}')ds' \quad (6)$$

Where $\vec{J}(\vec{r}')$ and $\rho(\vec{r}')$ represent the current density and charge density on the surface respectively, $k$ is the wave-number, $\vec{r}$ and $\vec{r}'$ represent the observer and source location points, $\mu$ is the permeability and $\varepsilon$ the permittivity of the background material. The current density is modelled by Rao–Wilton–Glisson (RWG) basis function $f_s$ [12] and Galerkin testing strategy is employed. Applying the continuity equation (4) the resultant MoM matrix entry for the $i^{th}$ observation and $j^{th}$ source basis is given by:

$$Z(i,j) = \frac{j\omega\mu}{4\pi}\iint f_i \frac{e^{-jk|\vec{r}-\vec{r}'|}}{|\vec{r}-\vec{r}'|}.f_j\,ds'\,ds +$$
$$\frac{1}{j\omega 4\pi\varepsilon}\iint \nabla f_i \frac{e^{-jk|\vec{r}-\vec{r}'|}}{|\vec{r}-\vec{r}'|}.\nabla f_j\,ds'\,ds \quad (7)$$

The MoM matrix thus formed is dense and presents a time and memory bottleneck. An iterative solution of the MoM matrix can be accelerated exploiting the physics of Green's Function interactions using a variety of techniques. Methods like MLFMA [4, 5] are kernel dependent whereas low-rank matrix compression based methods [6-10] are purely algebraic. The memory and solve time complexity of MoM can be mitigated by adapting low-rank matrix compressed H-Matrix decomposition. For the H-Matrix construction, the compression scheme can be applied on an oct-tree based 3D geometry decomposition [37, 38], where the matrix compression is applied for block interaction satisfying the admissibility condition given below

$$\min(dia(\Omega_t), dia(\Omega_s)) \leq \eta\,dis(\Omega_t,\Omega_s) \quad (8)$$

The admissibility condition states the minimum of the block diameter of the test $(\Omega_t)$ and source block $(\Omega_s)$ should be less than or equal to the admissibility constant $(\eta)$ times the distance between the test and source blocks. The oct-tree partition of the geometry is carried out until the number of elements in the block is less than *leaf size* [39]. At the leaf level, the block interaction not satisfying admissibility condition is considered as a near-field interaction. In the case of the multilevel oct-tree, the far-field block interacted at a higher level will not interact at the lower level. Fig. 1(a) and 1(b) shows the oct-tree, level 3, H-Matrix partition for the two-dimensional strip and a square plate where the green blocks are the admissible blocks (far-field blocks) at different levels and can be compressed by using various compression technique like ACA. The red blocks are the non-admissible blocks (near-field blocks) at the leaf level and form the dense near-field matrix interaction.

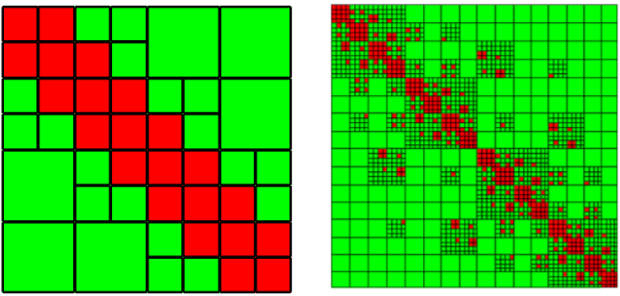

**Fig. 1**. *H-matrix partition of (a) Two-dimensional strip and (b) a square plate for oct-tree level 3*

In this work, the re-compressed ACA method [32, 33] is employed for the computation of H-Matrix. For the $m \times n$ rectangular sub-matrix $Z_{sub}^{m\times n}$ representing the coupling between two well-separated groups of $m$ observer bases and $n$ source bases, in the full MoM matrix, the ACA algorithm aims to approximate $Z_{sub}^{m\times n}$ by $A^{m\times n}$. In particular, the algorithm constructs $A^{m\times n}$ in product form as:

$$A^{m\times n} = U^{m\times r} \times V^{r\times n} \quad (9)$$

Where $r$ is the effective rank of the matrix $Z_{sub}^{m\times n}$ such that $r \ll \min(m,n)$, $U^{m\times r}$ and $V^{r\times n}$ are two dense rectangular matrices, such that:

$$\|R^{m\times n}\| = \|Z_{sub}^{m\times n} - A^{m\times n}\| \leq \varepsilon\|Z_{sub}^{m\times n}\| \quad (10)$$

For a given tolerance $\varepsilon$, where $R$ is termed as the error matrix and $\|.\|$ refers to the matrix Frobenius norm. Traditional ACA-based methods suffer from higher rank and error for the desired tolerance. To overcome this, a re-compression scheme is suggested in [31], where the orthogonal components in $U$ and $V$ are extracted and re-compressed by applying Singular-Value-Decomposition (SVD). The compression cost for each sub-matrix is given by $r^2(m+n)$ and the storage and matrix-vector product cost by $r(m+n)$. This compression scheme is applied in the present work on an oct-tree based geometry decomposition using a multilevel algorithm to reduce the overall matrix setup cost to $O(N\log N)$ and matrix-vector product cost to $(N\log N)\times p \times n_p$. However, the mere application of the re-compression ACA is not enough since for a large number of right-hand-side (RHS) vectors $n_p$, the algorithm efficiency can rapidly degrade if the number of iterations $p$ is high, necessitating an effective preconditioning strategy. Further processing is required as described in the following:

## 3. Schur's Complement Preconditioner

In this section, an efficient computation of diagonal block preconditioner from a symmetric near-field matrix using Schur's complement method is presented. Conventional H-Matrix near-field interaction is not symmetric in nature and symmetric property of the near-field matrix can be achieved by incorporating procedure given in [28]. The Schur's complement method [27] is a popular method in linear algebra for diagonalizing and solving block matrices. In [28], Schur's complement method is applied on symmetric block MoM matrices for decoupling the diagonal blocks in domain decomposition framework. In this work, Schur's complement method is applied to the near-field matrix, scaling it to diagonal block form. In the following subsection, the details of Schur's complement preconditioner computation process along with the ways to make it faster for setup and matrix vector product is presented.

### 3.1. Near-field Schur's Complement

For the preconditioner computation first, the geometry is divided into blocks based on the same oct-tree as used in the compression algorithm. In a multi-level H-Matrix compression scheme, as described in section 2, the MoM matrix $[Z]$ can be represented as a combination of the near-field interaction $[Z_N]$ at the leaf level and the compressed far-field interaction $[Z_F]$ obtained at multiple levels of far-interaction.

$$[Z]x = [Z_N + Z_F]x = b \quad (11)$$

Where $x$ represents the unknown coefficient vector, $b$ is the excitation vector. To explain the Schur's complement preconditioner setup, a demonstrative leaf-level cube structure comprised of 4 cubes, as shown in Fig. 2, is considered.

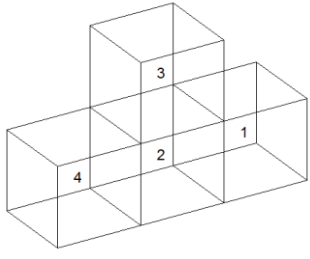

**Fig. 2**. *Representative leaf-level cubes for illustration of Schur's process*

For Fig. 2, block interaction between 1 and 4 forms the far-field interaction and rest form the near-field interaction. So near-field block matrix $[Z_N]$ in the case of Fig. 2 is given as:

$$[Z_N] = \begin{bmatrix} Z_{11} & Z_{12} & Z_{13} & 0 \\ Z_{21} & Z_{22} & Z_{23} & Z_{24} \\ Z_{31} & Z_{32} & Z_{33} & Z_{34} \\ 0 & Z_{42} & Z_{43} & Z_{44} \end{bmatrix} \quad (12)$$

Now the near-field can be scaled completely to a diagonal block by using the left and right scaling coefficients. Right scaling coefficient $[\alpha_1]$ for scaling the first row blocks $Z_{11}$ and $Z_{13}$ can be represented as:

$$[\alpha_1] = \begin{bmatrix} I_{11} & \alpha_{12} & \alpha_{13} & 0 \\ 0 & I_{22} & 0 & 0 \\ 0 & 0 & I_{33} & 0 \\ 0 & 0 & 0 & I_{44} \end{bmatrix} \quad (13)$$

Where, $I_{11}$ $I_{22}$ $I_{33}$ and $I_{44}$ are the identity block matrices. In the null-field method [25, 26] either left or right scaling is performed for the near-field scaling to diagonal blocks ignoring the fill-in blocks. In contrast, in the proposed Schur's complement method both left and right scaling are performed simultaneously considering fill-in blocks for complete near-field scaling to diagonal block. The fill-in blocks in the scaling matrix are further compressed and are discussed in further sub-sections. For row scaling $[\alpha_1]$ can be given as:

$$[\alpha_1] = \begin{bmatrix} I_{11} & -Z_{11}^{-1}Z_{12} & -Z_{11}^{-1}Z_{13} & 0 \\ 0 & I_{22} & 0 & 0 \\ 0 & 0 & I_{33} & 0 \\ 0 & 0 & 0 & I_{44} \end{bmatrix} \quad (14)$$

Similarly, for the complete scaling of column blocks $Z_{21}$ and $Z_{31}$ the left scaling coefficient $[\alpha'_1]$ is used and $[\alpha'_1]$ can be given as:

$$[\alpha'_1] = \begin{bmatrix} I_{11} & 0 & 0 & 0 \\ \alpha'_{12} & I_{22} & 0 & 0 \\ \alpha'_{13} & 0 & I_{33} & 0 \\ 0 & 0 & 0 & I_{44} \end{bmatrix} \quad (15)$$

$$[\alpha'_1] = \begin{bmatrix} I_{11} & 0 & 0 & 0 \\ -Z_{21}Z_{11}^{-1} & I_{22} & 0 & 0 \\ -Z_{31}Z_{11}^{-1} & 0 & I_{33} & 0 \\ 0 & 0 & 0 & I_{44} \end{bmatrix} \quad (16)$$

Now, equations (12), (13) and (16) can be combined to scale the first row and column block of $[Z_N]$ to diagonal block and the system of the equation can be given as:

$$[\widetilde{Z}_N^1] = \begin{bmatrix} I_{11} & 0 & 0 & 0 \\ \alpha'_{12} & I_{22} & 0 & 0 \\ \alpha'_{13} & 0 & I_{33} & 0 \\ 0 & 0 & 0 & I_{44} \end{bmatrix} \begin{bmatrix} Z_{11} & Z_{12} & Z_{13} & 0 \\ Z_{21} & Z_{22} & Z_{23} & Z_{24} \\ Z_{31} & Z_{32} & Z_{33} & Z_{34} \\ 0 & Z_{42} & Z_{43} & Z_{44} \end{bmatrix}$$
$$\times \begin{bmatrix} I_{11} & \alpha_{12} & \alpha_{13} & 0 \\ 0 & I_{22} & 0 & 0 \\ 0 & 0 & I_{33} & 0 \\ 0 & 0 & 0 & I_{44} \end{bmatrix} \quad (17)$$

Equation (17) can be represented as:

$$[\widetilde{Z}_N^1] = [\alpha'_1][Z_N][\alpha_1] \quad (18)$$

Performing the block multiplication in equation (17), $[\widetilde{Z}_N^1]$ can be represented as a block matrix form as:

$$[\widetilde{Z}_N^1] = \begin{bmatrix} Z_{11} & 0 & 0 & 0 \\ 0 & Z_{22}-Z_{21}Z_{11}^{-1}Z_{12} & Z_{23}-Z_{21}Z_{11}^{-1}Z_{13} & Z_{24} \\ 0 & Z_{32}-Z_{31}Z_{11}^{-1}Z_{12} & Z_{33}-Z_{31}Z_{11}^{-1}Z_{13} & Z_{34} \\ 0 & Z_{42} & Z_{43} & Z_{44} \end{bmatrix} \quad (19)$$

Equation (19) gives the Schur's complement of the first row and column block near-field matrix. Likewise, each row and column block can be scaled to form a diagonal block matrix and is of the form as given in equations below.

$$[\widetilde{Z}_N] = [\alpha'_3][\alpha'_2][\alpha'_1][Z_N][\alpha_1][\alpha_2][\alpha_3] \quad (20)$$

$$[\widetilde{Z}_N] = \begin{bmatrix} Z_{11} & 0 & 0 & 0 \\ 0 & \widetilde{Z}_{22} & 0 & 0 \\ 0 & 0 & \widetilde{Z}_{33} & 0 \\ 0 & 0 & 0 & \widetilde{Z}_{44} \end{bmatrix} \quad (21)$$

Equation (21) gives the complete scaled diagonal form of the near-field matrix. For solving the complete system of equations with left and right scaling coefficients the final system of equation (11) can be represented as:

$$[\alpha'_3][\alpha'_2][\alpha'_1][Z][\alpha_1][\alpha_2][\alpha_3][\widetilde{x}] = [\widetilde{b}] \quad (22)$$

Now, $[b]$ and $[x]$ in equation (11) can be extracted by

$$[\widetilde{b}] = [\alpha_3^T][\alpha_2^T][\alpha_1^T][b] \quad (23)$$

$$[x] = [\alpha_1][\alpha_2][\alpha_3][\widetilde{x}] \quad (24)$$

Equation (22) can be defined as the sum of the near and far-field as in equation (10) and is given as:

$$[\alpha'_3][\alpha'_2][\alpha'_1][Z_N + Z_F][\alpha_1][\alpha_2][\alpha_3][\widetilde{x}] = [\widetilde{b}] \quad (25)$$

Where, $[Z_F]$ is the far-field compressed ACA matrix blocks and $[Z_N]$ is the dense near-field block matrices. Equation (25) can be further simplified by as:

$$[\alpha'_3][\alpha'_2][\alpha'_1][Z_N][\alpha_1][\alpha_2][\alpha_3][\widetilde{x}] + \\ [\alpha'_3][\alpha'_2][\alpha'_1][Z_F][\alpha_1][\alpha_2][\alpha_3][\widetilde{x}] = [\widetilde{b}] \quad (26)$$

The first part of the above equation represents the block diagonal near-field as in equation (20) and then the equation can be further simplified as:

$$[\tilde{Z}_N][\tilde{x}] + [\alpha'_3][\alpha'_2][\alpha'_1][Z_F][\alpha_1][\alpha_2][\alpha_3][\tilde{x}] = [\tilde{b}] \quad (27)$$

Where, $[\tilde{Z}_N]$ is a scaled near-field diagonal block matrix which reduces the near-field matrix vector product time. However, due to the block diagonal nature of $[\tilde{Z}_N]$, it can be used as a preconditioner thus giving no fill-ins during factorization. The final preconditioned system of equation can be represented as

$$[\tilde{Z}_N]^{-1}\left[[\tilde{Z}_N] + [\alpha'_3][\alpha'_2][\alpha'_1][Z_F][\alpha_1][\alpha_2][\alpha_3]\right][\tilde{x}] = [\tilde{Z}_N]^{-1}[\tilde{b}] \quad (28)$$

An efficient preconditioner must have a low set-up, factorization and solve time. Diagonal block nature of the proposed preconditioner gives low factorization and solve time. In the next subsection, methods to reduce preconditioner set-up and overall matrix vector product time is presented.

### 3.2. Symmetric Schur's Complement

The preconditioner computation time can be reduced by exploiting the symmetric property of the near-field matrix. Due to the symmetric nature of the near-field matrix in equation (12) block matrix $Z_{21}$ will be transpose of $Z_{12}$ and $Z_{31}$ will be transpose of $Z_{13}$ i.e. lower diagonal blocks will be transpose of upper diagonal blocks. This gives a symmetric scaling coefficient blocks in the schur's preconditioner computation process. In the blocks scaling coefficients left scaling coefficient $[\alpha'_1]$ in equation (14) will be transpose of right scaling coefficient $|\alpha_1|$ in equation (16). Similarly $[\alpha'_2]$, $[\alpha'_3]$ and $[\alpha'_4]$ will be transpose of scaling coefficient $|\alpha_2|$, $|\alpha_3|$ and $[\alpha_4]$. Therefore, we can save scaling coefficient computation time by only computing right hand side coefficients. Time saving in precondition computation due to symmetric property of near-field matrix is shown in Table 2 for $5\lambda \times 5\lambda$ with 5985 unknowns.

**Table 2.** Precondition setup time

|  | Time (s) |
|---|---|
| Non Symmetric Near-field | 140.02 |
| Symmetric Near-field | 90.59 |

### 3.3. Near-field Block ordering

In the present approach, the generation of fill-in blocks may be efficiently reduced by using Graph ordering algorithms. Graph ordering algorithms are well known methods to reduce fill-ins in a sparse matrix [40]. Using this approach, the following example is presented.

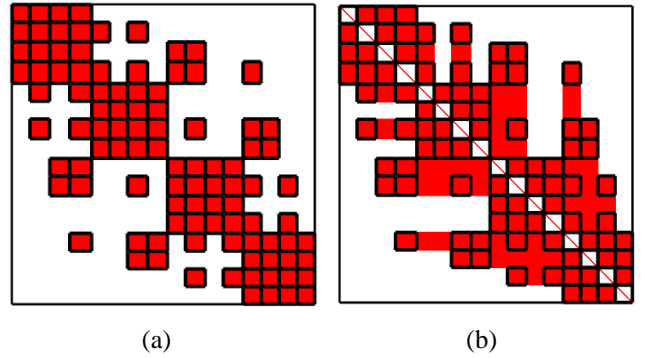

(a) (b)

**Fig. 3.** *(a) Near-field block matrix (b) Combined left and right scaling coefficients without ordering*

Consider the near-field block of a $5\lambda \times 5\lambda$ plate with 5985 unknowns, as shown in Fig. 3(a), computed from the leaf level oct-tree near-field bock interaction. By using left and right field scaling coefficients Fig. 3(b) can be generated. However, it can be further modified using the bandwidth reduction algorithm as described in Algorithm 1 to reduce precondition computation time, fill-in blocks and Number of Non Zeros (NNZ) elements in the scaling coefficients matrix.

**Algorithm 1**: Near-field Block Ordering

| | |
|---|---|
| *Step 1:* | *Assign Index to all the nodes at the leaf level Oct- tree which forms the nodes of the graph* |
| *Step 2:* | *Graph edges are filled based on the near-field node interaction at the leaf level* |
| *Step 3:* | *Graph algorithm for bandwidth reduction is applied* |
| *Step 4:* | *Basing upon then new graph node index is updated* |
| *Step 5:* | *Mesh element of triangle indices are updated for each node* |

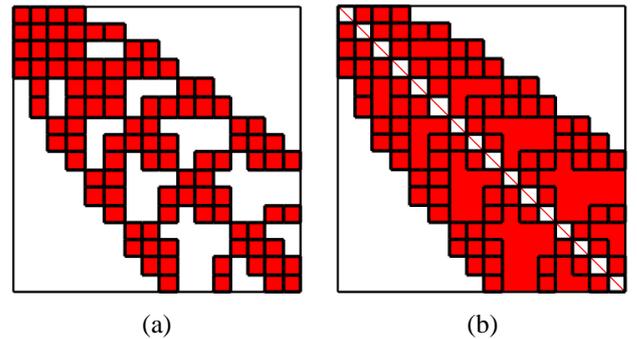

(a) (b)

**Fig. 4.** *(a) Near-field block matrix (b) Combined left and right scaling coefficients with Cuthill-Mikee ordering*

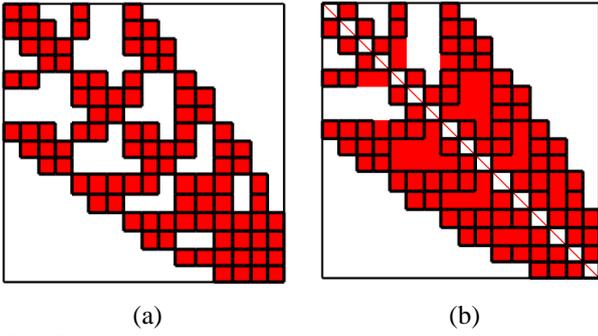

**Fig. 5.** *(a) Near-field block matrix (b) Combined left and right scaling coefficients with Reverse Cuthill-Mikee ordering*

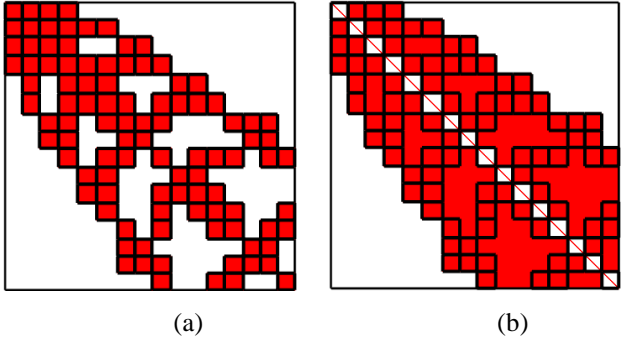

**Fig. 6**. *(a) Near-field block matrix (b) Combined left and right scaling coefficients with King's ordering*

The sparsity pattern generated by Cuthill-Mikee (Fig 4), Reverse Cuthill-Mikee (Fig 5), King's (Fig 6), and Sloan's (Fig 7) graph ordering algorithms for bandwidth reduction on the near-field and scaling coefficients blocks are shown. It can be observed from the scaling coefficient sparsity pattern figures that the ordering of near-field blocks affects the fill-in's in the scaling coefficient matrix.

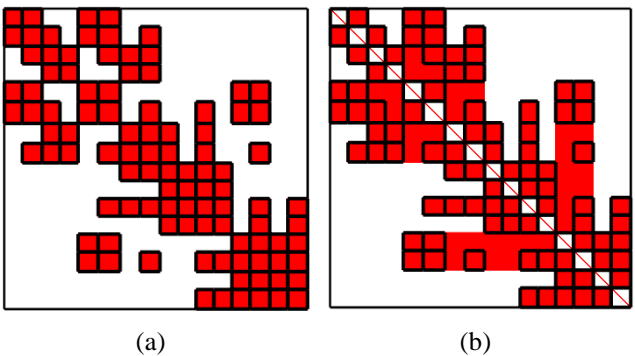

**Fig. 7.** *(a) Near-field block matrix (b) Combined left and right scaling coefficients with Sloan's ordering*

Effect of various graph ordering algorithm on the precondition computation time is shown in Table 3. Here all the precondition computation exploits the symmetric property of the near-field matrix. It can be observed that precondition computation time is reduced to half the convention computation time (Table 1) using Sloan's ordering algorithm.

**Table 3**. Precondition set-up time

| Graph Ordering | Time (s) |
| --- | --- |
| --- | 90.1414 |
| Cuthill-Mikee | 123.2073 |
| Reverse Cuthill-Mikee | 72.8518 |
| King's | 123.1237 |
| Sloan's | 71.3759 |

**Table 4**. Scaling Coefficient NNZ and Mat-Vec Time

| Graph Ordering | Oct-Tree(2) | | Oct-Tree(3) | |
| --- | --- | --- | --- | --- |
| | NNZ | Mat-Vec Time(s) | NNZ | Mat-Vec Time(s) |
| --- | 12458028 | 0.014804 | 5275335 | 0.011576 |
| Cuthill-Mikee | 12730712 | 0.015036 | 5576282 | 0.012470 |
| Reverse Cuthill-Mikee | 12346152 | 0.014649 | 5075171 | 0.011508 |
| King's | 12775196 | 0.016496 | 5710696 | 0.012697 |
| Sloan's | 12337473 | 0.014363 | 5144096 | 0.011240 |

Likewise, the effect of the graph ordering on near-field block and scaling coefficients are summarized in table 4 for oct-tree level 2 and 3. It can be observed that graph ordering of near-field blocks has a significant effect on precondition computation time and memory. Sloan's graph ordering gives better performance in terms of precondition computation time, memory and scaling coefficient matrix vector product time.

### 3.4. Compression of fill-in blocks in scaling coefficient matrix

The major cost of the iterative solver consists of the matrix vector product time. In the proposed method, the scaling coefficient vector product time can take a significant chunk of the overall matrix vector product time as shown in equation (26). It is well known that far-field block matrices can be compressed efficiently while maintaining the accuracy of the final solution. In this section, it is demonstrated that the scaling coefficient fill-in block matrices can also be similarly compressed using low-rank techniques for the given threshold. This leads to the reduction of scaling coefficient matrix vector time and storage. Considering the previous example, $5\lambda \times 5\lambda$ plate of 5985 unknowns with uniform meshing, the following observations are in order.

Left and right scaling coefficient combined sparsity pattern along with near-field blocks for $5\lambda \times 5\lambda$ plate are shown in Fig. 3(b) to 7(b). The red region outside the square box in Fig. 3(b) to 7(b) show the fill-in blocks in the scaling coefficients which can be compressed efficiently. Table 5 shows the compressed rank along with row and column size for the compression tolerance of 1e-3 and 1e-2 for some of the scaling coefficient fill-in blocks from Fig. 3(b).

**Table 5.** Compressed Scaling Coefficient Rank

| S. No | Row | Column | Rank(1e-3) | Rank(1e-2) |
|---|---|---|---|---|
| 1 | 363 | 385 | 18 | 13 |
| 2 | 385 | 385 | 20 | 12 |
| 3 | 363 | 363 | 17 | 13 |
| 4 | 385 | 363 | 21 | 14 |
| 5 | 363 | 363 | 23 | 13 |
| 6 | 363 | 363 | 16 | 11 |
| 7 | 385 | 363 | 20 | 13 |
| 8 | 363 | 363 | 21 | 14 |
| 9 | 363 | 385 | 19 | 12 |
| 10 | 363 | 385 | 17 | 11 |

Saving in terms of NNZ's in scaling coefficient and scaling coefficient matrix vector product time is shown in Table 6 for 1e-3 and 1e-2 compression tolerance for oct-tree level 2 and 3. It can be observed that compression of the fill-in blocks in scaling coefficients give memory savings as well as time saving in scaling coefficient matrix vector product time.

**Table 6.** Scaling Coefficient NNZ and Mat-Vec Time

| | Oct-Tree | NNZ | Mat-Vec Time(s) |
|---|---|---|---|
| No Compression | 2 | 16839988 | 0.019313 |
| Compression (1e-3) | 2 | 12458028 | 0.014804 |
| Compression (1e-2) | 2 | 12266298 | 0.014433 |
| No Compression | 3 | 8616302 | 0.016844 |
| Compression (1e-3) | 3 | 5275335 | 0.011576 |
| Compression (1e-2) | 3 | 4497586 | 0.009997 |

## 4. *Complexity Analysis*

In this section, the linear order complexity for preconditioner setup time, preconditioner matrix-vector product time and memory are demonstrated. For the complexity analysis, a uniform distribution of $N$ RWG bases in 3D is considered. The basis functions are grouped in a cube and following a multi-level oct-tree decomposition, each cube is recursively subdivided into 8 cubes starting from level 0 to level $L$. Therefore, at the lowest level there are $8^L$ leaf-level cubes. Assuming a uniform distribution, the number of basis functions in each leaf-level cube is $\frac{N}{8^L}$. Also, following the theory of most fast solver algorithms, it can be shown that for optimal efficiency of matrix storage and matrix-vector product cost $L = \log_8 N$. The preconditioner setup cost includes: (a) computation of left and right scaling coefficients $[\alpha']$ and $[\alpha]$ as in equation (14) and (16) which can be represented as $C_{SCC}$, (b) the second cost of preconditioner computation includes the scaling of the near-field to the diagonal block form by $[\alpha'][Z_N][\alpha]$ operation, for each row and column blocks in equation (18) and (20) which can be represented as $C_S$. Therefore, the total cost can be summed up as:

$$C_{TOTAL} = C_{SCC} + C_S \qquad (29)$$

### 4.1. Scaling coefficient computation cost

For the scaling coefficient computation, the major cost includes the inversion ($C_{MI}$) cost for diagonal block and the solving the inverse ($C_{SOL}$) for the row and column block near-fields as in equation (14) and (16). Therefor $C_{SCC}$ can be further be divided as the summation of inversion and solution cost as:

$$C_{SCC} = C_{MI} + C_{SOL} \qquad (30)$$

Inversion cost includes the single matrix inversion of a diagonal block for scaling near-field of each row and column blocks. Therefore, the matrix inversion cost of one matrix of a diagonal block at leaf level is given as:

$$C_{MI}^1 = k_1 \times \left[\frac{N}{8^L}\right]^3 \qquad (31)$$

Where, $k_1$ is a constant, the total cost for matrix inversion for the leaf level blocks is given by

$$C_{MI} = \sum_{i=1}^{8^L} C_{MI}^i = k_1 \times \left[\frac{N}{8^L}\right]^3 \times 8^L \qquad (32)$$

$$C_{MI} = k_1 \times N = O(N) \qquad (33)$$

For the computation of scaling coefficient, the inverted matrix has to be solved for the all row and column near-field blocks. Cost of matrix solution for one block at leaf level can be given as:

$$C_{SOL}^1 = k_2 \times \left[\frac{N}{8^L}\right]^2 \times \left[\frac{N}{8^L}\right] \qquad (34)$$

Where, $k_2$ is a constant, for a 3D structure, each block is surrounded by 26 near-field blocks. Therefore, the matrix solution cost for each row is given by:

$$C_{SOL}^{1R} = k_2 \times \left[\frac{N}{8^L}\right]^2 \times 26 \times \left[\frac{N}{8^L}\right] \qquad (35)$$

The total cost of the matrix solution at the leaf level blocks is the summation of the cost of each row and is given as:

$$C_{SOL} = \sum_{i=1}^{8^L} C_{SOL}^{iR} = k_2 \times \left[\frac{N}{8^L}\right]^2 \times 26 \times \left[\frac{N}{8^L}\right] \times 8^L \qquad (36)$$

$$C_{SOL} = k_2 \times 26 \times N = O(N) \qquad (37)$$

### 4.2. Near-field scaling cost

For the complete scaling of the near-field to diagonal block, the near-field blocks have to be multiplied by the left and right scaling coefficient by performing $[\alpha'][Z_N][\alpha]$ operation in equation (17) and (19) for each row. This procedure can be performed in two steps. First, computing $[Z_N][\alpha]$ and let it be represented by $[Z'_N]$ and the second process includes $[\alpha'][Z'_N]$ which leads to the complete diagonalization of the near-field matrix. Right scaling in equation (18) for each row can be given as

$$[Z'_1] = \begin{bmatrix} Z_{11} & Z_{12} & Z_{13} & 0 \\ Z_{21} & Z_{22} & Z_{23} & Z_{24} \\ Z_{31} & Z_{32} & Z_{33} & Z_{34} \\ 0 & Z_{42} & Z_{43} & Z_{44} \end{bmatrix} \begin{bmatrix} I_{11} & -Z_{11}^{-1}Z_{12} & -Z_{11}^{-1}Z_{13} & 0 \\ 0 & I_{22} & 0 & 0 \\ 0 & 0 & I_{33} & 0 \\ 0 & 0 & 0 & I_{44} \end{bmatrix} \quad (38)$$

$$[Z'_1] = \begin{bmatrix} Z_{11} & 0 & 0 & 0 \\ Z_{21} & Z_{22}-Z_{21}Z_{11}^{-1}Z_{12} & Z_{23}-Z_{21}Z_{11}^{-1}Z_{13} & Z_{24} \\ Z_{31} & Z_{32}-Z_{31}Z_{11}^{-1}Z_{12} & Z_{33}-Z_{31}Z_{11}^{-1}Z_{13} & Z_{34} \\ 0 & Z_{42} & Z_{43} & Z_{44} \end{bmatrix} \quad (39)$$

Now for the first block diagonalization left scaling has to be performed. Now the equation (17) can be written as:

$$[\widetilde{Z}_1] = \begin{bmatrix} I_{11} & 0 & 0 & 0 \\ -Z_{21}Z_{11}^{-1} & I_{22} & 0 & 0 \\ -Z_{31}Z_{11}^{-1} & 0 & I_{33} & 0 \\ 0 & 0 & 0 & I_{44} \end{bmatrix} [Z'_1] \quad (40)$$

Equation (39) decouples the $Z_{11}$ from the rest of the near-field block thus leading to a diagonal block formation. The major cost of the diagonalization comes from the right-hand scaling equation (38) which leads to the Schur's complement equation (39). Equation (38) mainly performs block matrix mutilation of the near-field of the first block with the right scaling coefficient. Therefore, the first-row scaling cost is given as

$$c_S^1 = k_3 \times \left[26\frac{N}{8^L} \times 26\frac{N}{8^L}\right] \times \left[\frac{N}{8^L} \times 26\frac{N}{8^L}\right] \quad (41)$$

Where $k_3$ is a constant, therefore the total cost of scaling the near-field blocks with the right scaling coefficients blocks can be represented as

$$C_S = \sum_{i=1}^{8^L} C_S^i \approx O(N) \quad (42)$$

In the case of left scaling in equation (40) $Z_{21}$ and $Z_{31}$ can be replaced by null blocks in equation (38) thus saving the cost. The above three sub-sections theoretically demonstrate that the overall cost of preconditioner generation, for a uniform 3D distribution of basis functions is $O(N)$. This experimentally shown in Fig. 8 for the increasing number of unknowns and the size of the plate structure. Further, the cost can be reduced by exploiting the symmetric property of equation (39)

### 4.3. Scaling Coefficient Matrix-Vector Product Cost:

Each left and right scaling coefficient consists of mainly the near-field blocks and come fill-in blocks which is compressed for 1e-2 threshold. Therefore, the scaling coefficient matrix vector product cost is given by:

$$C_{SMVP} = k_5 \times \left[\frac{N}{8^L} \times 26\frac{N}{8^L}\right] \times 8^L \quad (43)$$

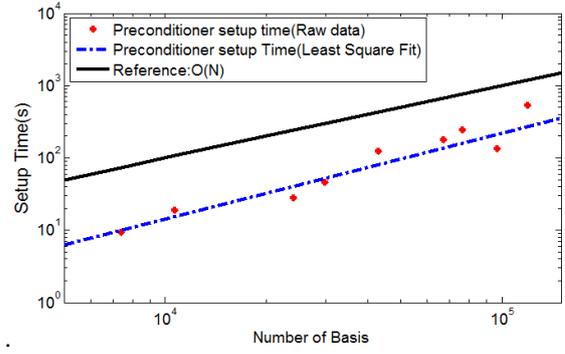

**Fig. 8.** *Preconditioner Setup Time*

Since the preconditioner matrix is a diagonal block in nature and the factorization and solve time of the diagonal block matrix in of $O(N)$. The $O(N)$ complexity of scaling coefficient matrix vector product time is shown experimentally in Fig.9 for the increasing number of unknowns and size of the plate structure.

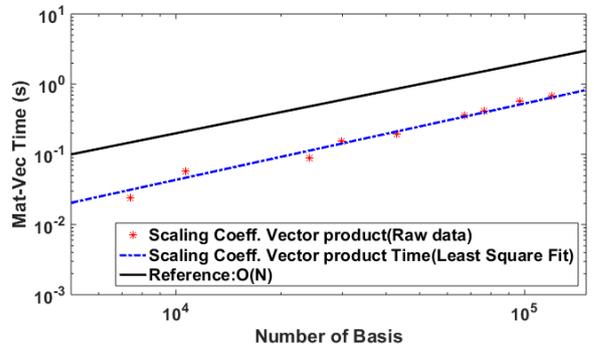

**Fig. 9.** *Left and right scaling coefficients combined matrix-vector product time*

### 4.4. Memory: Scaling coefficient and Preconditioner

The final matrices to be stored along with the compressed matrices includes the diagonal near-field which is used as a preconditioner matrix and the group of left and right -scaling coefficient matrix which sparse in nature. Each row of the scaling coefficient matrix consists of 26 sub-matrix blocks and 1 identity block matrix all of the size $N$. Therefore, the storage cost for one row-block of the scaling coefficient matrix is given by:

$$M_{SC}^1 = k_4 \times \left[\left[26\frac{N}{8^L} \times \frac{N}{8^L}\right] + \frac{N}{8^L}\right] \quad (44)$$

Where $k_4$ is the storage coonstant for each matrix elements. The total storage for the scaling coefficient can be represented as:

$$M_{SC} = \sum_{i=1}^{8^L} M_{SC}^i = k_4 \times \frac{N}{8^L}\left[\frac{26N}{8^L} + 1\right] \times 8^L \quad (45)$$

$$M_{SC} = k_4 \times 27 \times N = O(N) \quad (46)$$

The storage of block diagonal preconditioner is of $O(N)$. The $O(N)$ memory complexity of scaling coefficient matrix is shown experimentally in Fig.10 for the increasing number of unknowns and size of the plate structure.

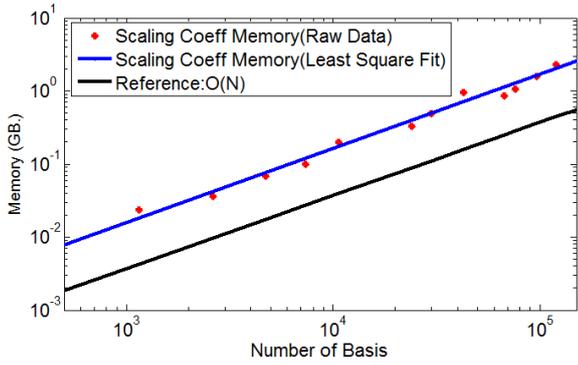

**Fig.10.** *Left and right scaling coefficients combined memory*

## 5. Numerical Results

In this section, numerical results are presented to validate the accuracy and efficiency of the proposed method. For all case studies, a GMRES tolerance of 1e-6 is considered. In the case of Schur's complement computation, Sloan's graph ordering is used to order the near-field blocks and the compression tolerance in scaling coefficients is kept as 1e-2. All the simulations for the results presents in this work are carried for double precision data type on the system with 128GB memory and Intel Xeon E5-2670 processor.

### 5.1. Validation of the proposed method with commercial tool

In this example, the proposed preconditioner is applied in alliance with an in-house re-compressed ACA-based H-Matrix. The results are validated with a commercial solver [41].

#### 5.1.1. Microstrip Line

A shorted microstrip with air background is considered with trace-width 0.05mm, length 10.0 mm and thickness 15.0µm. The ground plane is of width 1.0 mm, length 10.0 mm and thickness 15µm. Trace to ground plane distance is 30.0µm. The structure is meshed using uniform triangular elements generating 9681 RWG basis functions. In Fig.11 and 12 the magnitude and phase of S11 and S12 parameters are correlated with those generated from a commercial tool [40] using direct solver for EFIE dense matrix in a broadband simulation from 1 to 20 GHz.

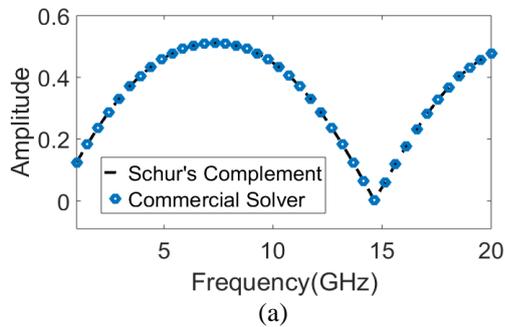

(a)

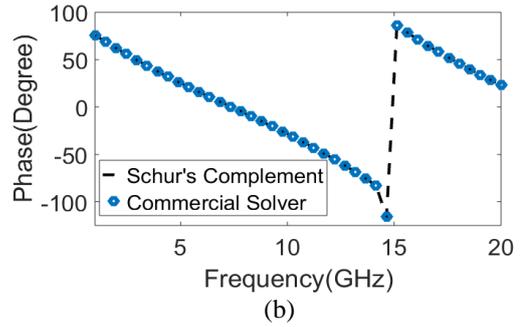

(b)

**Fig. 11**. *(a) Magnitude and (b) Phase of S11 for Microstrip Line*

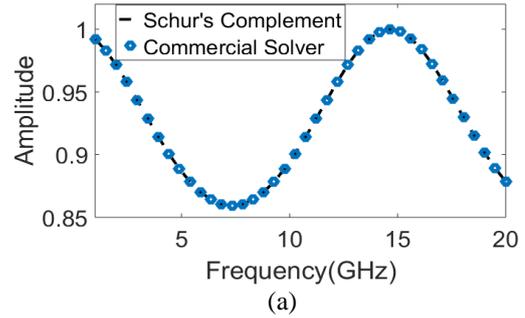

(a)

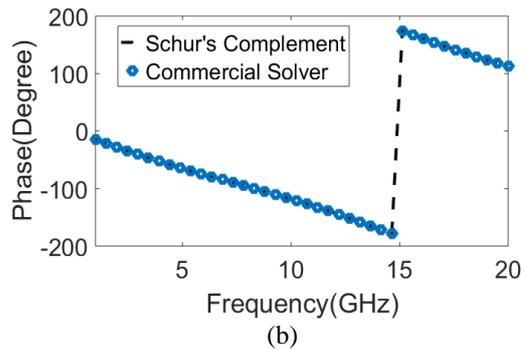

(b)

**Fig. 12**. *(a) Magnitude and (b) Phase of S12 for Microstrip Line*

The Eigen-value distribution before and after near-field Schur's complement preconditioning is shown in Fig. 13 for the MoM matrix at 1GHz.

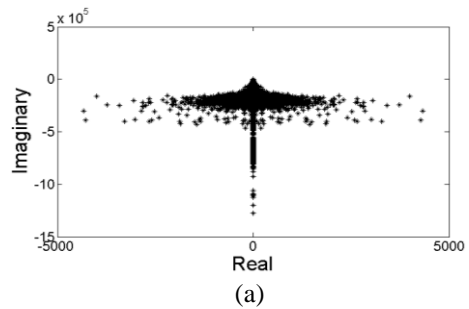

(a)

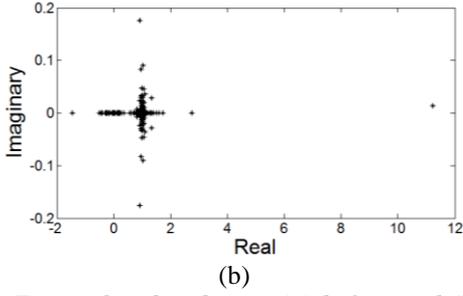

**Fig.13.** *Eigenvalue distribution (a) before and (b) after preconditioning of MoM matrix*

### 5.1.2. Aircraft

In this example, a model aircraft of length 4m and wingspan of 5m is taken. The structure is meshed using uniform triangular elements of size $0.1\lambda$ generating 83307 RWG basis functions. In Fig. 14 the monostatic RCS at 1GHz for plane wave incident at $\theta = 90°$ and $\phi = 0°$ to $180°$ is shown. The numerical result from our method compares well with the result obtained from the commercial solver [40]. In the case of commercial solver, EFIE MLFMA was used for matrix filling and for solving Bi-CGSTAB with stopping tolerance of 1e-6 was used. For the commercial solver along with ILU preconditioner, the total time taken to solve 180 RHS is 1150 hour.

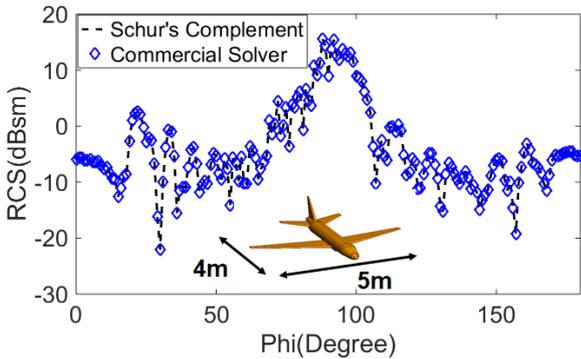

**Fig. 14**. *Monostatic RCS of Model aircraft with VV polarized plane wave incident with $\theta = 90°$ and $\phi = 0°$ to $180°$ at 1 GHz.*

### 5.2. Comparison with ILUT [23] and Null-Field Preconditioner [26]

In this section, the efficiency of the proposed method is validated. The performance of the proposed preconditioner is compared with that of ILUT and the null-field method. For ILUT the parameters are chosen as given in [24] and null-field with far-field assist [26] is used for the comparison. As discussed in [26] the relative efficiency of a preconditioner depends on some key parameters: (a) $t_{sm}$: MoM matrix setup time, (b) $t_{sp}$: preconditioner setup time, (c) $p$: average number of iterations required for convergence for 1 RHS, (d) $n$: number of RHS vectors, (e) $t_{mm}$: MoM matrix-vector product time and (f) $t_{mp}$: preconditioner operation time per iteration, which for the Schur's complement preconditioner includes $t_{mpp}$: preconditioner solve time and $t_{mps}$: scaling coefficient matrix-vector product time. The total matrix-setup and solve cost is given by:

$$t_{\text{total}} = t_{sm} + t_{sp} + [p \times n \times (t_{mm} + t_{mp})] \quad (46)$$

### 5.2.1. Microstrip Meander Line (MML)

A shorted meander line microstrip with air background, as shown in Fig. 15, is considered with the trace of width 0.05mm, length 10.0 mm and thickness 15.0µm. The ground plane is of width 5.0 mm, length 10.0 mm and thickness 15µm. Trace to ground plane distance is 30.0µm. The trace length is kept as 4.0mm and each bend-width is kept as 0.5mm. The speed-up comparison is done at 5 GHz and 40 GHz with different mesh size and frequency on the same structure. The comparison is shown in Table 7.

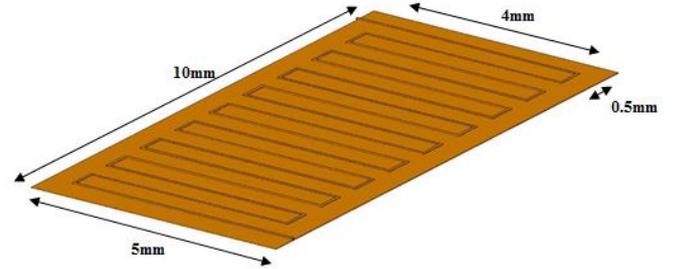

**Fig. 15.** *Microstrip meander line structure*

### 5.2.2. Multiport Traces (MT)

In this subsection, a multiple port via problem, as shown in Fig. 16, is considered. The height and diameter of via are 1.0mm and 1.0µm respectively and the antipad size is 0.3mm in diameter. The trace-width is 50µm and the thickness of the trace and ground plane layers are 20µm. The size of the ground plane is 10mm× 5.0 mm and the distance between the ground planes is 50µm for the set of 2 ports. The comparison is done for 4 ports at 20GHz and is shown in Table 7.

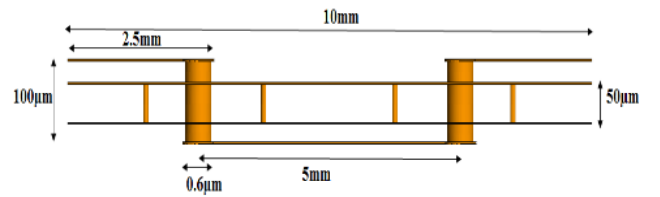

**Fig.16.** *Cross section view of multiport traces transitioning through via.*

### 5.2.3. Aircraft (AC)

In this subsection a model aircraft with dimensions given in section 5.1.2 is considered for the speed up comparison meshed at 1 GHz with $0.1\lambda$ mesh size. The comparison is shown in Table 7.

### 5.2.4. Cube (CU)

In this subsection a cube of length $6.0\lambda$ is considered for the speed up comparison meshed at 1 GHz with $0.1\lambda$ mesh size. The comparison is shown in Table 7.

**Table 7.** Schur's Complement Preconditioner Efficiency Comparison For different Mesh Sizes

| Problem | N | Preconditioner | $t_{sp}$ (s) | $p$ | $t_{mp}$ | | $t_{mm}$ (s) | $p \times (t_{mm} + t_{mp})$ (s) | $t_{total}$ (h) (All Ports / RHS) | Speed-up of Schur's PC (All Ports) |
|---------|---|----------------|--------------|-----|----------|---|--------------|----------------------------------|-----------------------------------|-------------------------------------|
|         |   |                |              |     | $t_{mpp}$ (s) | $t_{mps}$ (s) |          |                                  |                                   |                                     |
| MML (40GHz) 2 Ports | 109557 | Schur's    | 7402.886  | 138   | 0.144839 | 2.307487 | 0.65329  | 427.1961    | 3.87    | ----   |
|                     |        | Null-Field | 11659.149 | 358   | 0.151119 | 2.589863 | 2.58986  | 1908.4414   | 6.11    | 1.5 x  |
|                     |        | ILUT       | 16915.437 | 16691 | 3.881351 | ----     | 2.44588  | 105607.8126 | 65.18   | 16.8 x |
| MT (20GHz) 4 Ports  | 35685  | Schur's    | 847.594   | 57    | 0.105656 | 0.731477 | 0.209189 | 59.6400     | 0.726   | ----   |
|                     |        | Null-Field | 3544.430  | 223   | 0.105867 | 0.688537 | 0.538892 | 297.3250    | 1.739   | 2.39 x |
|                     |        | ILUT       | 2066.052  | 5431  | 0.395828 | ----     | 0.539810 | 5081.4499   | 6.943   | 9.5 x  |
| AC (1GHz) 180 RHS   | 83307  | Schur's    | 1453.143  | 1580  | 0.079789 | 0.619148 | 0.838416 | 2429.0177   | 124.50  | ----   |
|                     |        | Null-Field | 2708.948  | 2338  | 0.079810 | 0.655274 | 1.137914 | 4379.0693   | 222.36  | 1.78 x |
|                     |        | ILUT       | 6951.729  | 11539 | 0.901223 | ----     | 1.130787 | 23447.3633  | 1176.95 | 9.45 x |
| CU (300MHz) 180 RHS | 87303  | Schur's    | 19513.004 | 524   | 0.046505 | 1.097207 | 0.880525 | 1060.7001   | 82.98   | ----   |
|                     |        | Null-Field | 25416.310 | 676   | 0.047350 | 1.348644 | 1.680345 | 2079.6051   | 135.56  | 1.63 x |
|                     |        | ILUT       | 25515.164 | 10123 | 1.601163 | ----     | 1.686106 | 33277.0240  | 1695.46 | 20.43 x |

## 6. Conclusion

In this paper, a new preconditioning technique is proposed based on the Schur's complement method. The proposed precondition computation method scales complete near-field in a compact block-diagonal form with an associated left and right scaling matrices. The scaled block-diagonal matrix is further used as a preconditioner. 2x speed for precondition computation is achieved by exploiting the symmetric property of near-field and using graph ordering algorithm. The saving will be further exaggerated with an increase in the number of RHS vectors. The over all cost of the solve time is dominated by the preconditioner set up time which can be resolved by efficient parallelization which will further improve the speed-up.